# An E-PINN assisted practical uncertainty quantification for inverse problems

Xinchao Jiang[1,a], Xin Wang[a], Ziming Wen[a], ENying Li[c], Hu Wang[2,a,b]

a. *State Key Laboratory of Advanced Design and Manufacturing for Vehicle Body, Hunan University, Changsha, 410082, P.R. China*

b. *Beijing Institute of Technology Shenzhen Automotive Research Institute, Shenzhen 518000, People's Republic of China*

c. *College of Mechanical & Electrical Engineering, Central South University of Forestry and Technology, Changsha, 410082, P.R. China*

## Abstract

How to solve inverse problems is the challenge of many engineering and industrial applications. Recently, physics-informed neural networks (PINNs) have emerged as a powerful approach to solve inverse problems efficiently. However, it is difficult for PINNs to quantify the uncertainty of results. Therefore, this study proposed ensemble PINNs (E-PINNs) to handle this issue. The E-PINN uses ensemble statistics of several basic models to provide uncertainty quantifications for the inverse solution based on the PINN framework, and it is employed to solve the inverse problems in which the unknown quantity is propagated through partial differential equations (PDEs), especially the identification of the unknown field (*e.g.*, space function) of a given physical system. Compared with other data-driven approaches, the suggested method is more than straightforward to implement, and also obtains high-quality uncertainty estimates of the quantity of interest (QoI) without significantly increasing the complexity of the algorithm. This work discusses the good properties of ensemble learning in field inversion and uncertainty quantification. The effectiveness of the proposed method is demonstrated through several numerical experiments. To enhance the robustness of models, adversarial training (AT) is applied. Furthermore, an adaptive active sampling (AS) strategy based on the uncertainty estimates from E-PINNs is also proposed to improve the accuracy of material field inversion problems.

---

[1] First author. *E-mail address*: jiangxinchao@hnu.edu.cn (X.C. Jiang)

[2] Corresponding author. Tel: +86 0731 88655012; fax: +86 071 88822051.
   *E-mail address*: wanghu@hnu.edu.cn (H. Wang)

## 1. Introduction

The inverse problem is an important area of science and industry. Its main purpose is to identify the quantity of interest (QoI) through indirect observations or measurements inside the domain. Practically, due to the complexity of the given forward problem, the number of observations in an inverse problem might be limited. Furthermore, any sensor measurements are inevitably disrupted by noise, and the numerical round-off of the value might also affect the accuracy of results. These characteristics typically lead to inverse problems being ill-posed. It includes solution existence, solution uniqueness, and instability of the solution process [1]. All of these issues pose abundant challenges to inverse problems.

Over the past few decades, scholars have proposed numerous kinds of methods and strategies to handle inverse problems. Classical inverse problem research has focused on creating conditions to ensure the existence of solutions to such ill-posed problems, as well as a more stable approximation of solutions in the presence of noise [2, 3]. For example, applying some regularization methods to impose well-posedness, and then a relatively definite problem can be solved by optimization or other manners [2]. This kind of solution might not be the best answer to the original inverse problem. Therefore, it is necessary to address the inverse problem within the framework of statistical inference [4, 5], as uncertainty analysis of this approximated solution is usually required. A widely used method is the Bayesian inference [6, 7], which provides a robust way of taking parameter fluctuations into account. This framework formulates a complete probabilistic description of the unknown QoI and uncertainties given observed data. The posterior distribution of the QoI is obtained by the fusion of likelihood and prior information. However, obtaining an accurate posterior is often computationally expensive since the exploration of posterior distributions is generally based on sampling methods, such as Markov chain Monte Carlo (MCMC) [8, 9] or Sequential Monte Carlo (SMC) [10]. Admittedly, these classical methods have achieved great success in many problems, but whether optimization-based approaches or Bayesian inference may require numerous repeated simulations of the forward problem. Therefore, it is natural that the reductions or surrogates for the forward operator are employed to accelerate the computationally intensive problems [11-14]. However, such attempts at acceleration are still based on pre-prepared repeated forward simulations. It is computationally prohibitive in situations where the physical model is complex or the model spaces increase to high dimensions, and care should also be taken to avoid the so-called inverse crime [15]. There is no doubt that classical methods represent a major step forward in the field of inverse problems. But these methods are limited by the inability to capture bespoke structures (*e.g.*, underlaying physics) in data that vary in different applications [16]. For example, a typical requirement in the area of engineering and industry is the reliable recovery of hidden high-dimensional model parameters (*e.g.*, sources, materials, boundaries etc.) from noisy data. It makes classical methods difficult to implement in these situations.

With the rapid growth of data size and advances in computational science. Machine learning-based (ML) data-driven approaches have become promising methods to solve inverse problems, which show the impact on a wide range of inverse problems, such as medicine [17], heat transfer [18], and geophysics [19]. The data-driven approach concerned in this paper is deep neural networks (DNNs) which is a branch of ML. However, there is generally insufficient data to support purely data-driven methods in almost all engineering applications [20], which severely limits the implementation of purely data-driven approaches to many inverse problems. Whereupon a recent line of development in computational science has emerged and it has gradually become a hot research field -- physics-informed machine learning [21], especially the physics-informed neural networks (PINNs) in reference [22]. This novel computational paradigm shows how DNNs can solve both the forward and inverse problems in a unified framework. It empirically demonstrates how to learn the constraints of physical information from data and built reliable surrogates from incomplete parameters and relatively small data sets using modern computing tools [23, 24]. The underlying physical law of observations is coded into the loss of DNNs to impose regularization on the parameters of this physics- and data-driven surrogate model, meaning the bespoke structure of the data can be captured by PINNs. It sheds new light on solving inverse problems and naturally introduces regularization into data-driven approaches. The PINN is more than a promising way to solve high-dimensional inverse problems, it also avoids repeated simulations of the forward problem. Therefore, this study focuses on the application of PINNs in solving inverse problems. Lu *et al.* [25] proposed an approach that PINNs with hard constraints (hPINNs) to solve high dimensional PDE-constrained inverse problems (geometry design). The proposed hPINNs adopt an idea similar to that of full-space methods for PDE-constrained optimization, where the objective is optimized and the PDE is solved simultaneously. Yin *et al.* [26] employed PINNs successfully to extract the permeability and viscoelastic modulus form thrombus deformation data. Their work demonstrated that PINNs can infer the properties of material from noisy synthetic data. Haghighat *et al.* [27] introduced PINNs to learning and discovery in solid mechanics. Then the application of this method in parameter identification of linear elasticity and nonlinear elastoplastic problems is illustrated. Tartakovsky *et al.* [28] used PINNs to estimate the hydraulic conductivity of saturated and unsaturated flows governed by Darcy's law. Interestingly, the idea of "ensemble" emerged in their study. In their work, an ensemble of DNNs is trained to analyze the average relative error in response fields. These representative efforts show the wide application of PINNs to inverse problems. Considering the propagation of the model and data uncertainties, many scholars also extended the work of uncertainty quantification to the PINN. These extensions mainly include generative adversarial networks (GANs) [29, 30], Bayesian deep learning (BDL) [31], and polynomial chaos extension (PCE) [32]. These well-designed algorithms achieve satisfactory results in the studied problems, whereas they introduce greater computational complexity and are more difficult to implement. Compared with the original PINN, they are harder to train and more assumptions and hyperparameters that might limit the expressiveness of the model are introduced. Such as the computational cost of training GANs is much higher than training feed-forward neural networks and prone to mode-collapse [30]. BDL might face an intractable likelihood function, especially when there have no direct observations of the QoI. Reference [31] shows that the performance of the variational inference-based BDL and MC-dropout [33] is poor. Although the HMC-based BDL has achieved good results, the HMC algorithm is too sensitive to the number of steps and step size [8]. It is difficult to obtain appropriate settings for many problems. Hence, these issues raise the need for a simple and practical uncertainty quantification method to solve inverse problems under the framework of PINNs. Recently, Yang *et al.* [34] have proposed multi-output physics-informed neural networks (MO-PINNs) to obtain predictive uncertainty by modifying the output layer of PINNs. The MO-PINN is relatively simple to implement and the examples demonstrate the feasibility of this method. Nonetheless, the inverse problem presented in their study is only limited to the problem of parameter identification in low-dimensional space. Interestingly, MO-PINNs might also be interpreted as an ensemble where the predictions are averaged over an ensemble of neural networks with the same parameters except the output layer. Therefore, it motivates the exploration of the ensemble as an alternative solution for uncertainty quantification in PINNs. Moreover, extensive experience has shown that the ensemble of models can improve the predictive performance [35].

To address the above issues, this study introduces deep ensemble (DE) [36] into the framework of PINNs and proposes a simple and practical method named ensemble PINNs (E-PINNs) to solve the inverse problem with uncertainty quantification. The final predictions and uncertainty estimates are achieved through an ensemble of PINNs (different parameters). Notably, this is not the first-time ensemble learning has been explicitly used to improve the performance of PINNs. Reference [37] discussed how to improve the accuracy by using the ensemble method in training. However, the task of this study is to estimate the field rather than a few parameters, which usually increases the ill-posedness of the inverse problem as it attempts to recover an infinite-dimensional QoI from some indirectly noisy observations. Therefore, the significant difference is that our work focuses more on the reliability of ensemble learning in field inversion and discusses the good properties of ensemble learning in uncertainty quantification. In addition, the suggested method obtains good inversion results of unknown fields with relatively low computational complexity. No sophisticated tools are necessary, and it also avoids the problem of likelihood functions and uncertainty propagation which are difficult to handle in the Bayesian framework. The feasibility of E-PINNs is demonstrated by solving inverse problems where QoI is the unknown function (*i.e.*, spatial or temporal fields) and quantifying the unknown which is propagated through PDEs. The results are compared with MC-dropout since it is also relatively simple to implement. Furthermore, an adaptive active sampling (AS) strategy is proposed to improve the accuracy of results in material field inversion problems.

This paper is organized as follows. In Section 2, the formulation of the inverse problem is illustrated. In Section 3, the methodology is described and the details of E-PINNs are given. Several numerical examples are carried out in Section 4. The results demonstrate that the proposed method is effective in the source field and material field inversion. Finally, Section 5 summarizes the core of the present work. All code and data accompanying this study will be made publicly available at https://github.com/yoton12138 after the paper is officially published.

## 2. Problem Setup

Scholars frequently wish to relate physical parameters $\boldsymbol{m} \in M$ that characterize a model to the collected observation data $\boldsymbol{d} \in D$. It is usually assumed that the underlying physics associated with $\boldsymbol{m}$ and $\boldsymbol{d}$ are adequately understood. Hence, a function/operator $\mathcal{G}(\cdot)$ can be used to describe this relationship:

$$\mathcal{G}(\boldsymbol{m}) = \boldsymbol{d} . \qquad (1)$$

Here, $M$ (model parameter space) and $D$ (data space) are vector spaces with appropriate topologies, and their elements represent possible model parameters and data, respectively. Generally, the forward problem is to obtain *d* with given *m*, or $\mathcal{G}: M \to D$. Whereas this paper focuses on the inverse problem of identifying *m* with given *d*, or $\mathcal{G}^{-1}: D \to M$. A significant issue is that actual measurements always contain a certain amount of noise. The noise can come from sensors or any numerical round-off. Equation (1) usually takes the form with a noise component $\xi$:

$$d = d_{true} + \xi = \mathcal{G}(m_{true}) + \xi, \qquad (2)$$

where the forward operator and the probability distribution of the noise are derived from the first principles. $\mathcal{G}(\cdot)$ refers to an operator when *m* and *d* are functions and a function when *m* and *d* are vectors. Operators can be in many forms, such as ordinary differential equations (ODEs), partial differential equations (PDEs), or even a linear or nonlinear system of algebraic equations. It should be noted that this study considers the inverse problem in which *m* is a function, which needs to be distinguished from the parameter identification problem in which *m* is a low-dimensional vector or a certain value. Therefore, a physical system governed by PDEs defined on a domain $D \in R^d$ can also be described as:

$$\mathcal{O}_x(d(x); m) = 0, x \in D, \qquad (3)$$

with boundary conditions (BCs):

$$\mathcal{B}_x(d(x)) = 0, x \in \partial D. \qquad (4)$$

Here, $\mathcal{O}_x$ represents a set of operators, $\mathcal{B}_x$ is the general form of BC operators, and $\partial D$ is the boundary of the domain *D*. *m* is the parameter field that specifies the PDEs, which is the QoI of the inverse problem. The general purpose of this study is to solve the inverse problem concerning *m* by minimizing the objective function $\mathcal{L}$. From a modeling viewpoint, it can also be interpreted as constructing a surrogate model of *m* in the space-time domain by using finite noisy data not directly related to *m*.

## 3. Methodology

3.1 *The physics-informed neural networks used in this study*

Physics-informed neural networks (PINNs) are first proposed in the study of Raissi *et al.* [22], which demonstrate how classical physical laws enhanced the robustness and efficiency of modern machine learning algorithms. It introduces a new scientific modeling paradigm through a combination of physics-driven and data-driven. Different from the architecture of original PINNs, the schematic of the PINN employed in this study is shown in Figure 1. The input feature in a PINN is spatiotemporal variables in Cartesian coordinates, i.e., $\boldsymbol{x} = (x, y, t)$. In this study, to highlight the application of PINNs in field inversion, two fully connected neural networks $\hat{d}(\boldsymbol{x};\theta)$ and $\hat{m}(\boldsymbol{x};\vartheta)$ parameterized by $\theta$ and $\vartheta$ are used to model the observation data $\boldsymbol{d}$ (or measurements) and the QoI $\boldsymbol{m}$, respectively. In this paradigm, the process of solving the continuous inverse problem is approximated by finding a set of parameters $\vartheta$ to describe the QoI $\boldsymbol{m}$ best under given observation data. It is essentially a reparameterization of QoI. Assuming that each neural network has $L$ layers, the transformation of each layer can be written as follows:

$$\mathbf{h}^l = \mathcal{A}^l(\mathbf{W}^l \cdot \mathbf{h}^{l-1} + \mathbf{b}^l), l = 1, 2, \cdots, L, \qquad (5)$$

where $\mathbf{h}^0$ is the input layer and $\mathbf{h}^L$ is the output layer, $\mathbf{W}^l$ and $\mathbf{b}^l$ are weights and bias of each layer $l$. $\mathcal{A}^l$ is the activation functions, and default to be set to hyperbolic tangent (tanh) since it is smooth enough to compute arbitrary-order derivatives. Hence, such an explicit transformation is easy to obtain the gradients of the outputs concerning the inputs through automatic differentiation (AD) (*e.g.* Pytorch, TensorFlow). The restrictions on $\hat{d}(\boldsymbol{x};\theta)$ taking two steps in PINNs. Firstly, the surrogate model $\hat{d}(\boldsymbol{x};\theta)$ should be able to approximate the true observation data $\boldsymbol{d}$. Then along with a set of collocation points $\boldsymbol{x}_f$ (for PDE constraint) and $\boldsymbol{x}_{bc}$ (for BCs constraint) inside the spatiotemporal domain, one can restrict the two networks $\hat{d}(\boldsymbol{x};\theta)$ and $\hat{m}(\boldsymbol{x};\vartheta)$ to satisfy the PDEs and BCs, *i.e.*, $\mathcal{O}_x(d(\boldsymbol{x});\boldsymbol{m}) = 0$ and $\mathcal{B}_x(d(\boldsymbol{x})) = 0$. Hence, the general loss function of this problem can be written as:

$$\mathcal{L} = \omega_1 \mathcal{L}_{data} + \omega_2 \mathcal{L}_{pde} + \omega_3 \mathcal{L}_{bc}, \tag{6}$$

$$\mathcal{L}_{data} = \frac{1}{N_d} \left\| \hat{d}(\boldsymbol{x};\theta) - \boldsymbol{d} \right\|_2^2, \tag{7}$$

$$\mathcal{L}_{pde} = \frac{1}{N_f} \left\| \mathcal{O}_x(\hat{d}(\boldsymbol{x}_f;\theta); \hat{m}(\boldsymbol{x}_f;\vartheta)) \right\|_2^2, \tag{8}$$

$$\mathcal{L}_{bc} = \frac{1}{N_{bc}} \left\| \mathcal{B}_x(\hat{d}(\boldsymbol{x}_{bc};\theta)) \right\|_2^2, \tag{9}$$

where $N_d$ is the number of $\boldsymbol{d}$, $N_f$ and $N_{bc}$ are the number of $\boldsymbol{x}_f$ and $\boldsymbol{x}_{bc}$ respectively. $\omega = \{\omega_1, \omega_2, \omega_3\}$ is the weight of each loss term. The loss $\mathcal{L}$ can generally be optimized by the stochastic gradient descent (SGD) method. Here, adaptive moment estimation (Adam) [38] is adopted. As for collocation points, there are many ways to obtain them. Such as Latin hypercube sampling [39], Sobel sequence [40], or grid-based approaches [41]. To cover the pre-defined spatiotemporal domain more evenly, grid-based approaches are applied to sample these collocation points in this study.

*3.2 E-PINNs: deep ensemble-based physics-informed neural networks*

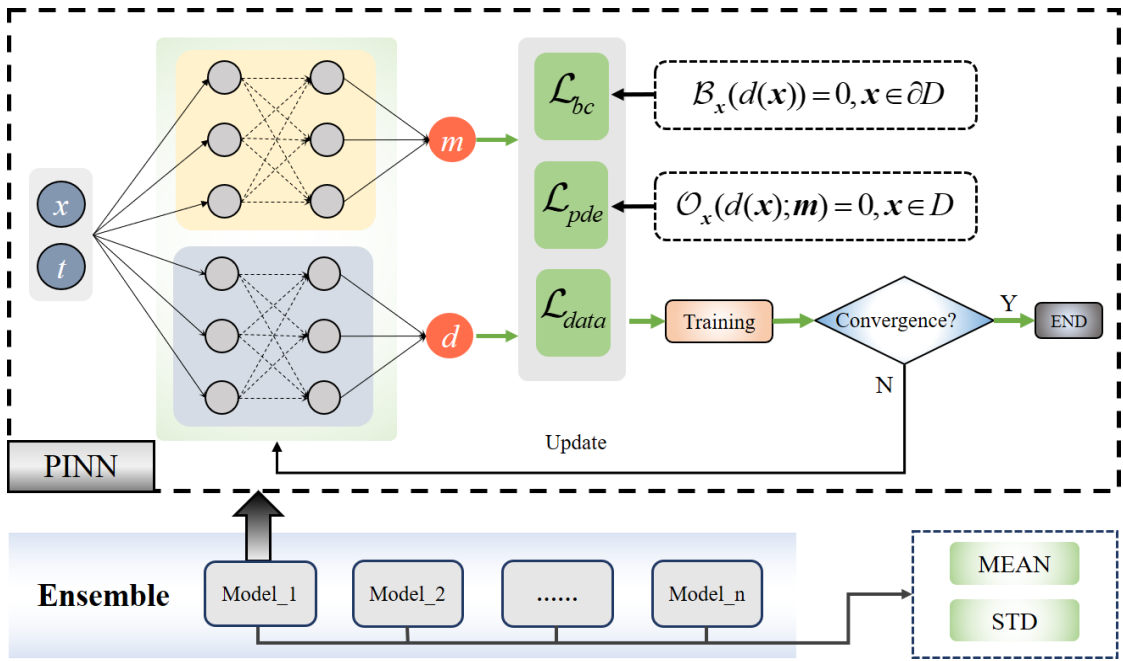

**Figure 1.** The schematic of the E-PINN for field inverse problems

DE is a powerful uncertainty quantification method, which has been widely used in many practical applications due to its simple implementation and scalability [42]. As previously discussed, exact Bayesian inference is often computationally intractable for DNNs because the distributions and likelihood functions of the datasets are unknown. Furthermore, a complex algorithm design will lead to implementation difficulties. Therefore, the DE naturally becomes an attractive method for uncertainty quantification. As shown in Figure 1, the DE is introduced to PINNs, and E-PINNs are proposed. The critical issues of the proposed method are adversarial training (AT) and ensemble learning, respectively.

*3.2.1 Adversarial training: fast gradient method*

AT [43-45] is an important training method in the field of deep learning, which is widely used in computing vision and natural language processing fields. It is a regularization technology, which improves the robustness of the model by attacking the model with adversarial samples. AT leads to a better generalization than conventional training. It is also a form of data augmentation, where an extra dataset $\mathcal{D}' = \{\mathbf{x} + \Delta\mathbf{x}, \boldsymbol{d}\}$ is obtained by applying a small perturbation to the original training dataset $\mathcal{D} = \{\mathbf{x}, \boldsymbol{d}\}$. The instances in $\mathcal{D}'$ are called adversarial samples. Here, an AT approach named the fast gradient method (FGM) [46] is introduced into the model training of PINNs. It is different from the fast gradient sign method (FGSM) [45] which uses sign functions to impose the same size perturbation on each dimension in the original DE. FGM scales each input dimension with a specific gradient and yields better adversarial samples than FGSM. In general, AT [47] can be written in the following format:

$$\min_{\theta} \mathbb{E}_{\theta(\boldsymbol{x},d) \sim \mathcal{D}}[\max_{\Delta\boldsymbol{x}} \mathcal{L}(\boldsymbol{x} + \Delta\boldsymbol{x}, d; \theta)] \tag{10}$$

where $\Delta\boldsymbol{x}$ is the perturbation, $\mathcal{L}(\boldsymbol{x} + \Delta\boldsymbol{x}, d; \theta)$ is the loss of each sample. In FGM, the way to compute the $\Delta\boldsymbol{x}$ in the procedure of maximizing the $\mathcal{L}$ is:

$$\Delta\boldsymbol{x} = \tau \frac{\nabla_x \mathcal{L}(\boldsymbol{x}, d; \theta)}{\|\nabla_x \mathcal{L}(\boldsymbol{x}, d; \theta)\|_2}, \tag{11}$$

where $\tau$ is a given small value such that the max-norm of the perturbation is bounded.

Hence, the total loss in the PINN will plus a new term with the same weight as $\mathcal{L}_{data}$:

$$\mathcal{L}_{AT} = \frac{1}{N_d}\left\|\hat{d}(\boldsymbol{x}+\Delta\boldsymbol{x};\theta)-\boldsymbol{d}\right\|_2^2. \tag{12}$$

The procedure of AT can be interpreted as a computationally efficient solution to smooth the predictive distribution via increasing the likelihood of the target around a small neighborhood of the original training samples [36]. AT eradicates the vulnerability in a single model by forcing it to learn robust features. To some extent, it can avoid overfitting in the context of a few noisy data.

### 3.2.2 Ensemble: training and prediction

In this study, E-PINNs are proposed to quantify the uncertainty in inverse problems. The solutions and uncertainty estimates are obtained through an ensemble of PINNs. Specifically, there are two ensemble schemes: randomization-based and boosting-based. The randomization-based scheme is employed in the proposed approach. Compared with the boosting-based approach, the randomization-based is easier to implement, and the members of this ensemble scheme can be trained in parallel without any interaction. It means E-PINNs can easily be used for distributed computing. Each PINN member of the E-PINN is trained using the entire dataset since DNNs typically perform better with more data. This has also been demonstrated in the study [48] that a complete dataset and randomly initialized neural networks can obtain better results than bagging and other methods for DE. Intuitively, using more data to train the base learners help reduce their bias. This also makes it easier to establish a physically consistent model under the PDEs constraints. It is even more significant in an inverse problem since more observations might avoid the dilemma of multi-solutions or ensure the solution uniqueness. The framework of E-PINNs is summarized below:

---

**Algorithm1** E-PINNs for field inverse problem with uncertainty quantification

**Step 1**: Given the training set $\mathcal{D}$ and perturbation $\tau$ (1% of the smallest input ranges), the number of models *M*, and other the hyperparameters of the E-PINN.

**Step 2**: Initialize parameters $\theta_1, \cdots \theta_M$ and $\vartheta_1, \cdots \vartheta_M$ randomly to construct the E-PINN.

**Step 3**: Specify a total loss function $\mathcal{L}$.

**Step 4**: Train each PINN to find the best parameters $\theta_i$ and $\vartheta_i$ independently (parallel is recommended), $i = \{1, \cdots, M\}$.
**Step 5**: Ensemble and make predictions of the QoI ***m***.

In this ensemble method, each member is evenly-weighted. The final inversion of the QoI of this study is the posterior $p(\boldsymbol{m}|\boldsymbol{x},\boldsymbol{d})$, which is considered a mixture of Gaussian distribution. Once the E-PINN is trained, meaning a surrogate model concerning the QoI is built. The solution of inverse problems is equivalent to the prediction of the model in such a framework. The prediction of the QoI is only related to the inputs. Hence, it can be approximated by $\hat{m}(\boldsymbol{x};\vartheta)$:

$$p(\boldsymbol{m} \mid \boldsymbol{x}) \sim \mathcal{N}(m(\boldsymbol{x}), \sigma^2(\boldsymbol{x})), \tag{13}$$

$$\mu(\boldsymbol{x}) \approx M^{-1} \sum \hat{m}(\boldsymbol{x};\vartheta_i), \tag{14}$$

$$\sigma^2(\boldsymbol{x}) \approx M^{-1} \sum \hat{m}^2(\boldsymbol{x};\vartheta_i) - \mu^2(\boldsymbol{x}). \tag{15}$$

By using ensemble statistics, the equation (14) is an unbiased estimation of the mean of QoI, and the equation (15) yields a maximum likelihood estimation (MLE) of the variance of a Gaussian posterior. The feasibility of this estimation method can be understood from the training mechanism of DNNs. Methods using an ensemble of models to produce predictive distributions rely on the model diversity to produce accurate uncertainty estimates [49]. The training objectives in this multi-task learning problem are diversative, where the surrogates need to satisfy with data distribution, physics information, and adversarial samples. The random initialization and SGD algorithms introduce the uncertainty of model parameters into the model training. This is particularly important when there is insufficient data to construct an accurate and robust model. Since different strategies (initialization, etc.) will lead to different approximations of local features of the unknown quantities. Therefore, this kind of diversified loss function and stochastic training method can obtain models with enough diversity, they tend to converge to different local optimal solutions. This is a natural way of considering the uncertainty of the model itself, also known as epistemic uncertainty. In addition, noisy data introduce an unquantifiable aleatory uncertainty into the suggested method. Although this study does not directly quantify the aleatory uncertainty, the use of AT allows noises to indirectly participate in the model training. The proposed approach also avoids dealing with the likelihood function and parameter adjustment which are difficult to handle in other algorithms (*e.g.,* BDL).

3.3 *Adaptive active sampling*

AS is an algorithm that guides the data sampling process. In the inversion process, especially in the ill-posed inverse problem, the solutions of QoI are often uncertain. A better result and the reduction of uncertainty can be achieved by AS at important locations. It is of great importance when data acquisition costs are high. In E-PINNs, it can be implemented by adding a new observation to the training of the model. AS is carried out in two steps. The first is to find the place the point $x_s$ with the largest variance $\sigma^2(x)$ in the current inverse solution. Then, sensors or other measurement manners are used to obtain a new precise measurement value $m_s$ of the point directly, and the value is added to the training set to update the parameters of $\hat{m}(x;\vartheta)$. Notably, the stop criteria for AS are worth discussing, because it is difficult to give convergence criteria for unknown QoI without any reference. The traditional method of using iteration error as the stopping criterion might lead to the problem of excessive sampling. For example, it is much more difficult to directly measure anisotropic material fields than to obtain response fields, the cost of AS needs to be controlled. Therefore, an adaptive stop criterion is proposed that the AS stops while the maximum uncertainty is less than a threshold value $\eta=\alpha\Delta$, $\alpha$ is a scaling factor which is suggested to be 0.02, $\Delta$ is the value range of $\mu(x)$. These settings consider a trade-off between accuracy and cost. The adaptive AS can be summarized as follow:

---
**Algorithm2** Adaptive active sampling

**Step 1**: Strat with a preliminary trained E-PINN with $\theta_1,\cdots\theta_M$, $\vartheta_1,\cdots\vartheta_M$, and $\eta$.

**Step 2**: **While** the maximum uncertainty $< \eta$ **do**

    Compute the standard deviation using equation (15) and find the location $x_s$.

    Obtain the observation $m_s$ of the QoI and add it to the training set.

    Update the parameters $\theta_1,\cdots\theta_M$ and $\vartheta_1,\cdots\vartheta_M$.
---

# 4. Numerical results and discussion

In this section, several numerical examples are presented to demonstrate the feasibility and effectiveness of E-PINNs. It should be noted that the purpose of this study is not to compare the computational cost with other algorithms, but to an alternative way to solve the ill-posed inverse problem and obtain reasonable uncertainty estimates simply and practically. The inverse problem involved in this section is crucial for many applications. The QoI is a spatial field, which may correspond to inhomogeneous material properties or may represent distributed source terms in govern equations.

4.1 *Source field inversion problem*

In this numerical example, a source inversion problem is used to demonstrate the effectiveness of E-PINNs. The source here can be considered to be contaminant sources or heat sources. The purpose is to identify the source field based on finite noisy observations obtained from a single forward simulation. Unlike many problems of estimating the partial coordinate or scale parameters, this is a spatial function identification problem. It greatly increases the ill-posedness of the inverse problem. Specifically, a diffusion system in a 2-dimensional spatial domain with 2 sources can be described as:

$$\lambda(\frac{\partial^2 u}{\partial x^2} + \frac{\partial^2 u}{\partial y^2}) + s = 0,$$
$$s = \sum_{i=1}^{2} k_i \exp(-0.5 \times \sigma_i^{-2} \times \|x_i - p_i\|^2), \quad (16)$$

where $u$ is the concentration of the source, namely the observation data, the diffusion coefficient $\lambda = 0.02$, and the $s$ is the unknown source, namely QoI, which needs to be identified. $k = \{1, 2\}$ and $\sigma = \{0.15, 0.05\}$ are the strength and scale parameters of each Gaussian source, respectively. The exact positions for the two sources are: $p = (p_1, p_2)$, $p_1 = (0.3, 0.4)$, $p_2 = (0.8, 0.8)$. The exact numerical solution solved by FEM with 4,960 triangle meshes and 2,484 nodes is shown in Figure 2. The outer boundary of this example is the Dirichlet boundary condition with a value of zero, and the inner boundary is adiabatic in the normal direction. Assuming there are some measurements collected through the sensors located at each node, the noise level of $u$ is $\xi_u \sim \mathcal{N}(0, 0.02^2)$. A subset of the measurements which contains 2,000 random samples for $u$ and 2,400 uniformly grid-based collocations points are utilized as the training data.

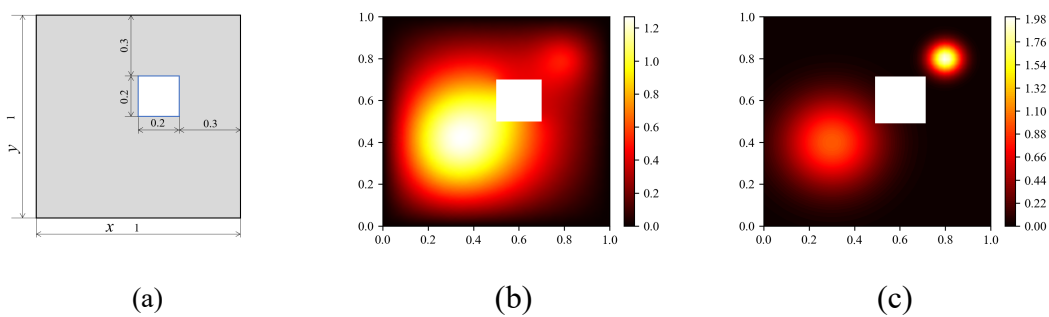

(a)                (b)                (c)

**Figure 2.** Basic information of source field inversion problem. (a) geometry description; (b) the exact numerical solution for the diffusion system; (c) the source field.

In this section, two fully connected neural networks with 4 hidden layers are used for each member of the E-PINN, each hidden layer contains 20 units. One for approximating the observed responses $u$ and the other for the QoI $s$. Note that there are no observations and prior information of the QoI, and its surrogate model is constructed by forcing DNNs to satisfy certain physical constraints. The hyperparameters of the Adam optimizer are $\beta_1 = 0.9$, $\beta_2 = 0.999$, and the different learning rates of $\hat{d}(\boldsymbol{x}; \vartheta)$ and $\hat{m}(\boldsymbol{x}; \vartheta)$ are 0.003 and 0.005, respectively. The total number of training steps of each model in the E-PINN is 100,000. In this example, the weight values in the loss function are all set to 1. As for MC-dropout, the same network architecture and activation function is used as the member of E-PINNs, and the dropout rate of 0.1 is added after each non-linearity in the surrogate model of QoI. These architectures of DNNs and other hyperparameters are determined empirically to yield accurate results. The optimal choice of these settings is another area of study and is outside the scope of this work.

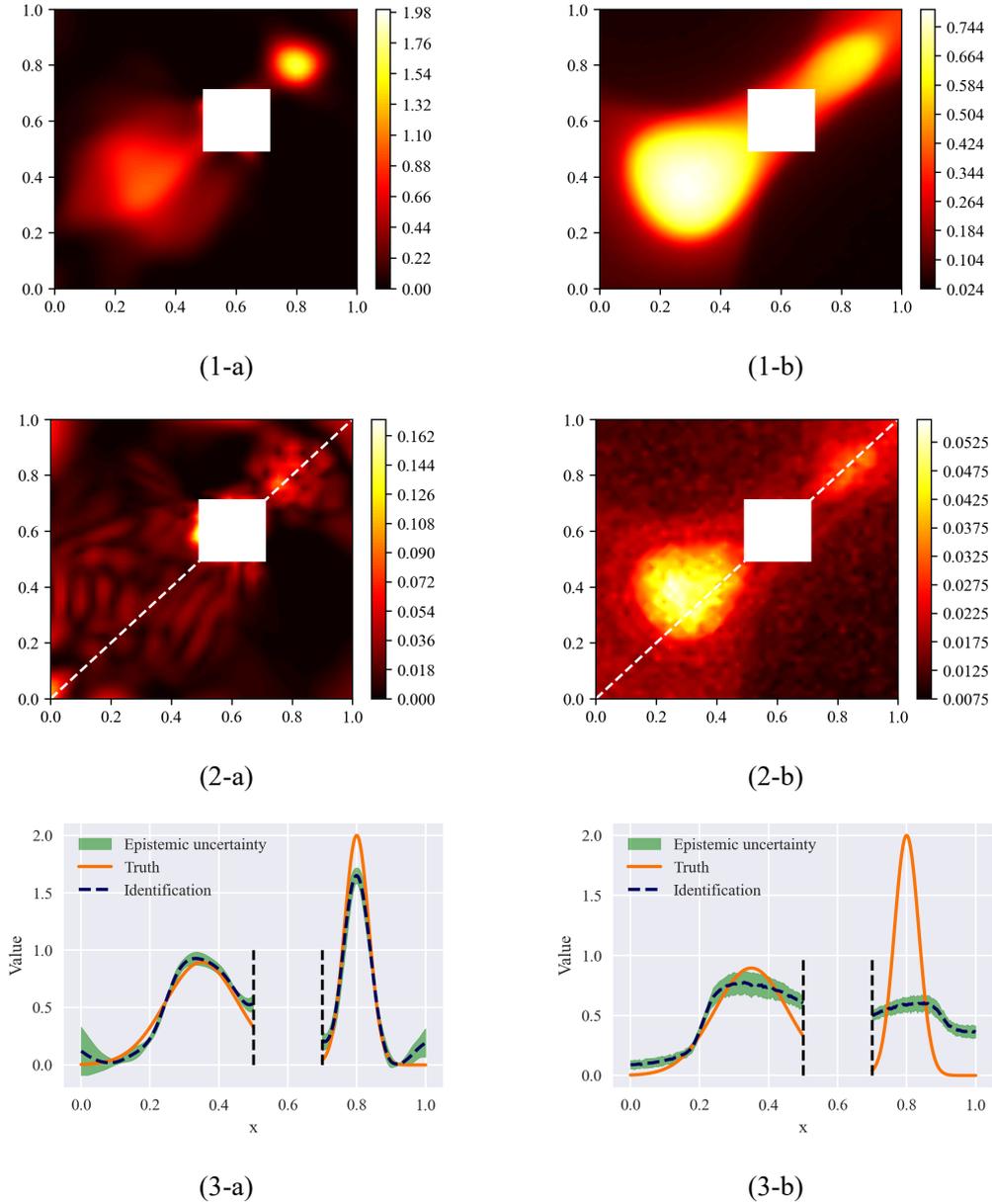

**Figure 3.** Identification of the source field. The first row is the results of the E-PINN, the second row is the results of MC-dropout; Column (a) is the mean of inversion; Column (b) is the standard deviation of inversion; Column (c) is the details of the diagonal domain.

The results are shown in Figure 3. As can be seen from the results, the E-PINN has an excellent ability to identify source fields and yields high-quality uncertainty estimates, which reaches a relative $L_2$ error of 0.1770 and an $R^2$ of 0.9568. The results here are obtained using the default value $M = 5$ in E-PINNs. However, it can be observed that the results of MC-dropout are far from the ground truth, which only reaches a relative $L_2$ error of 0.8646, and an $R^2$ of 0.7008. The E-PINN captures the main features of the source field, with errors mainly distributed at peak points and geometric boundaries. The epistemic uncertainty interval (twice the standard deviation, $\pm 2\sigma$) of the inversion is mainly concentrated around the sources and the corners of the boundaries since the evolution behavior of the response field in these areas is more complex. As can be seen from the diagonal detail plot, real values are mostly in the credibility interval, and the quality of uncertainty estimates is also better than that of MC-dropout. It should be mentioned that for noisy data, perfect inversion of the source field is challenging, which may require a much larger amount of data, and even more excellent algorithm design. After all, it is difficult to obtain the exact solution without any prior. For many problems, it is sufficient to identify the main characteristics of the source such as peak values, peak coordinates, and relatively accurate profile. Next, the effects of AT and the number of models in the E-PINN are studied, as shown in Figure 4 and Table 1. The results show that the E-PINN with $M = 5$ has sufficient accuracy and uncertainty quantification capability of the inverse solution. At least, simply model increasing does not improve the performance of E-PINNs in any meaningful way in this case, the accuracy fluctuates within a very small range. In E-PINNs, AT is used by default. According to the results of the controlled trials, the inverse results with AT perform better, but the improvement in this numerical example was not particularly obvious. In situations trained with AT, the inversion around the second source is more accurate and the peak value is closer to the real value. However, it also causes a certain reduction in the accuracy of the first source. An objective assessment of this result is that AT is a regularization training technique, like any regularization method, it does not play an obvious role in all problems. In the next example, the effect of AT will be shown more clearly. But in this case, it does not make the model much worse and therefore still worth considering as a more stable means to handle noisy data. This example also preliminarily verifies the feasibility of ensemble learning for PINNs to obtain credibility intervals by synthesizing the diversity of different models.

**Table 1.** Mean performance comparison of E-PINNs of different training schemes

| - | AT | | Without AT | |
|---|---|---|---|---|
| | $M=5$ | $M=10$ | $M=5$ | $M=10$ |
| $R^2$ | 0.9568 | 0.9534 | 0.9496 | 0.9492 |
| Relative $L_2$ error | 0.1770 | 0.1841 | 0.1912 | 0.1912 |

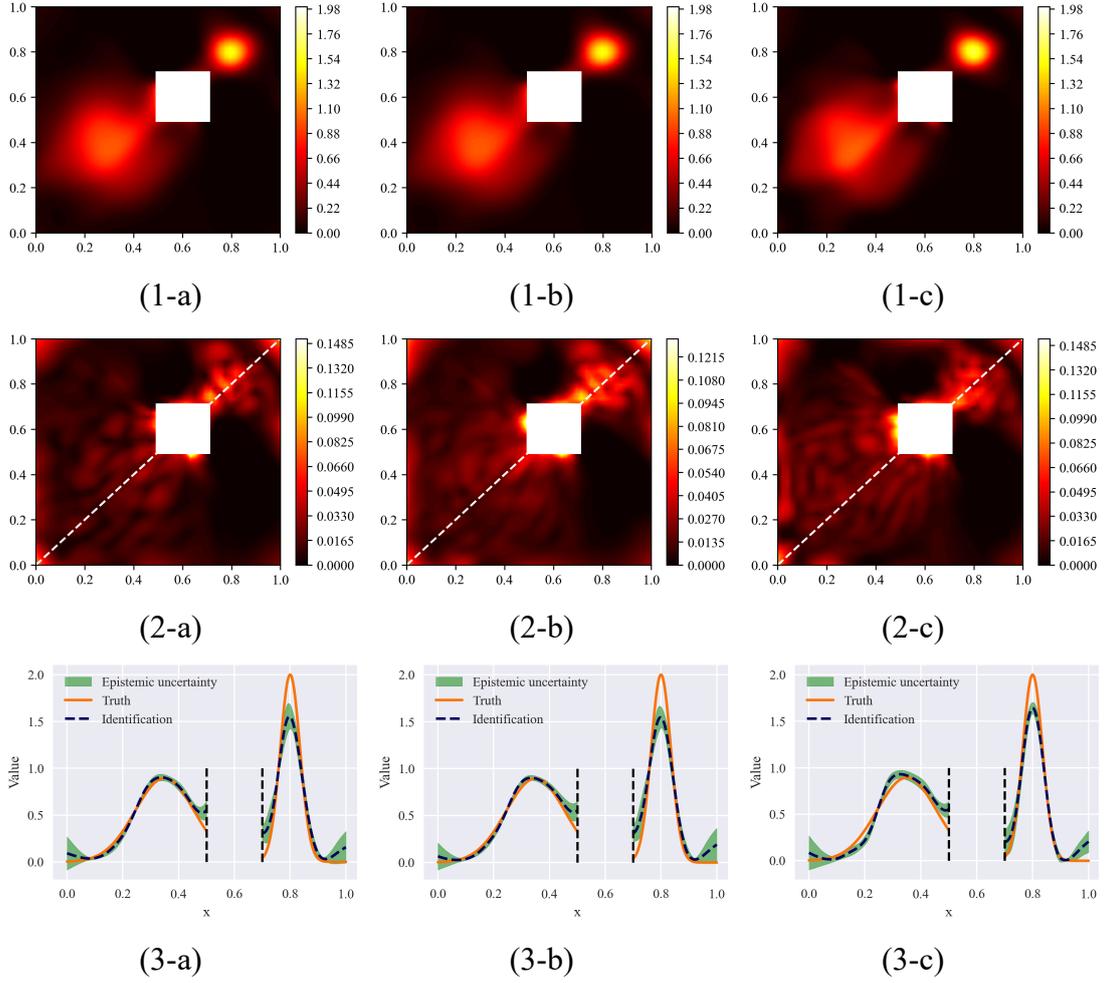

**Figure 4.** Controlled trials. The first row is the mean, the second row is the standard deviation, and the last row is the detail of the diagonal domain. Column (a) is trained without AT, $M=5$; Column (b) is trained without AT, $M=10$; Column (c) is trained with AT, $M=10$.

4.2 *Material field inversion in the nonlinear diffusion equation*

In this section, another example of material field inversion is considered to validate the accuracy and effectiveness of the proposed approach. This section is divided into two parts. The first part is to discuss the normal inversion procedure, and the second one is to study the adaptive AS method to improve accuracy.

4.2.1 *Inversion procedure*

This numerical example comes from reference [12]. The inverse problem setup is the same as the reference. In particular, a nonlinear transient diffusion equation on the interval $x \in [0,1]$ with adiabatic boundaries is considered:

$$\frac{\partial u}{\partial t} = \frac{\partial}{\partial x}(v(x)\frac{\partial u}{\partial x}) + \sum_{i=1}^{N=3} \frac{s_i}{\sqrt{2\pi}\sigma_i} \exp(-\frac{|l_i - x|^2}{2\sigma_i^2})[1 - H(t-0.01)],$$

$$\frac{\partial u}{\partial x}\bigg|_{x=0} = \frac{\partial u}{\partial x}\bigg|_{x=1} = 0, \qquad (17)$$

$$u(x, t=0) = 0.$$

The source term in equation (17) involves 3 sources, each source acts on the interval $t \in [0, 0.01]$ and is centered at $l_i \in \{0.25, 0.5, 0.75\}$ with a constant strength $s = 100$ and a variance $\sigma^2 = 10^{-3}$. $H(\cdot)$ is the Heaviside function. $v(x)$ is the material field, and supposing $v(x) > 0.1$. It can be considered as a prototype for the inverse estimate of an inhomogeneous material field, such as a thermal conductivity field or the permeability field in a porous medium. The space-time domain is uniformly dispersed into a $49 \times 301$ grid. Two formats of $v(x)$ are studied and the exact numerical solutions are solved by FEM. The details of the forward problems are shown in Figure 5. In this case, noting that $v(x)$ is the QoI $m$ of the equation (1), which is solved by E-PINNs. Considering the consistency with the description in the literature, the same logarithmic diffusivity $m(x) = \log(v(x) - 0.1)$ is used to describe parametric forms and draw the plots. One is a sinusoidal profile and the other is a random sample of the zero-mean Gaussian process with a covariance function $\exp(-0.5 \times |x_1 - x_2|^2 \times 0.3^{-2})$, or random field for short. Assuming there are virtual sensors located at some nodes, the measurements $\mathcal{D} = \{u_1^*(x_1, t_1), \cdots, u_n^*(x_n, t_n)\}$ take place over the time interval $t \in [0.01, 0.03]$. In this example, two different noise levels are studied, $\xi_1 \sim \mathcal{N}(0, 0.01^2)$, $\xi_2 \sim \mathcal{N}(0, 0.05^2)$. Specifically, the sensors are evenly distributed across the domain with a space interval of 1/12, and the time interval of each sensor is 0.25s. Hence, a dataset containing 117 noisy observations is collected for inversion.

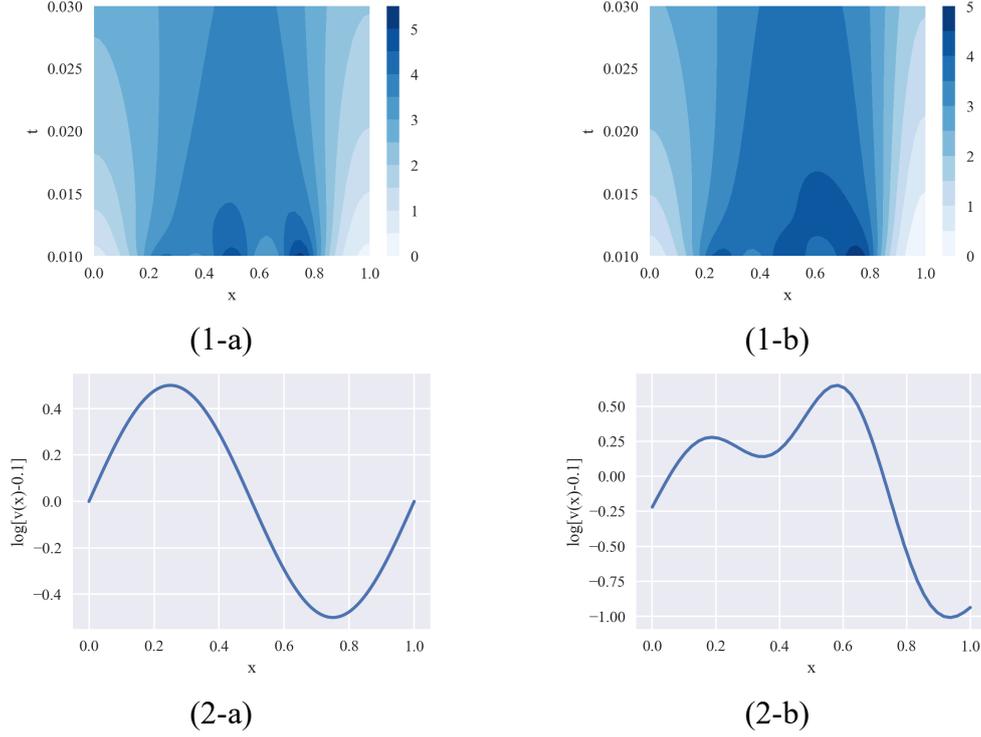

(2-a)　　　　　　　　　　　　　　　(2-b)

**Figure 5.** Numerical solutions of different $v(x)$. (1-a) $u$ field with the diffusivity in (2-a); (1-b) $u$ field with the diffusivity in (2-b).

Each model in the E-PINN has the same two neural networks as the previous example, one acting as a surrogate for response field $u$ and the other acting as a surrogate for QoI $v(x)$. From the comparisons of the previous example, the default value $M = 5$ can be used directly in this example. The number of training steps for each member of the E-PINN is 50,000. Other hyperparameters settings of Adam are the same as before. Since the observation set of this example is relatively small, a larger $50 \times 100$ uniformly grid-based collocation point set is used to strengthen the effect of the PDE constraint. $\omega_1$ is set to 1, $\omega_2$ and $\omega_3$ are set to 0.001. For MC-dropout, the same processing is maintained as before, with a fixed dropout rate of 0.1. It should be mentioned that some training tricks are used to accelerate convergence and improve accuracy in this example. Hierarchical training is introduced to accelerate the convergence of E-PINNs. In the framework of PINNs, once the inverse problem is solved, the surrogate of the forward operator is constructed exactly. Therefore, the surrogate of the forward operator $\hat{u}(x;\theta)$ could be trained with observations separately in the first step. A threshold $\delta = 10^{-4}$ is defined to ensure that each $\hat{u}(x;\theta)$ is sufficiently close to the distribution of observed data, meaning that the $\mathcal{L}_{data}$ of each surrogate is lower than $\delta$. This threshold can be equivalent to adaptive training steps to achieve expected convergence under different network weight initialization conditions. Then this pre-trained network is applied as the new "initialization" of a PINN. This is similar to an optimization procedure using a new and better initial iteration point. The second step is regular comprehensive multi-task training, including additional PDE and BC constraints.

Figure 6 and Table 2 show the inverse results of E-PINNs for different nonlinear diffusivity conditions and different noise levels after convergence. It can be seen that the mean and credibility interval of inverse results are reasonable. In all cases, the uncertainty of QoI is greater at the boundary, peak, and where the curvature is big, which is consistent with intuition. The results obtained here are not based on any assumptions about the diffusivity. However, it has achieved well performance, which reflects the efficiency and feasibility of E-PINNs. These results can be interpreted by the rich model diversity generated during the different training processes. It is difficult to construct an accurate and relatively stable model with limited information. The model constructed under different initializations tends to learn different local features of potential material fields that satisfy data distribution. Therefore, the ensemble model provides a more reasonable prediction, and the statistics of the E-PINN can be used to estimate the uncertainty. When the material field is sinusoidal, the overall accuracy is higher than that of the random field. This result may be the different resolution requirements of the observation set for different forms of nonlinear material fields. Figure 7 and Table 3 show the results of MC-dropout, these demonstrate that E-PINNs perform much better than MC-dropout in terms of all the performance measures. The obvious role of AT in this example is shown by comparing Figure 6 and Figure 8. In different cases, the accuracy and the quality of uncertainty estimates are improved to varying degrees. Because AT prevents the over-fitting of the model to noisy observation data and forces it to discover the more essential potential features in the data. For the sinusoidal profile, the $R^2$ is improved by about 1% to 2%. However, the credibility interval becomes narrower, which means the model is more confident in the inverse results. For the random field, the $R^2$ is improved by 5.1% and 7.0% when the noise is 0.01 and 0.05, respectively. The credibility interval is more consistent with the characteristics of real contour and the credibility interval contains more true values.

**Table 2.** The performance of E-PINNs in all cases

| Noise | Criterion | Sinusoidal profile | | Random field | |
|---|---|---|---|---|---|
| | | AT | Without AT | AT | Without AT |
| 0.01 | $R^2$ | 0.9911 | 0.9795 | 0.9336 | 0.8880 |
| | Relative $L_2$ error | 0.0945 | 0.1430 | 0.2576 | 0.3347 |
| 0.05 | $R^2$ | 0.9569 | 0.9377 | 0.9006 | 0.8415 |
| | Relative $L_2$ error | 0.2077 | 0.2496 | 0.3153 | 0.3981 |

**Table 3.** The performance of MC-dropout in all cases

| Noise | Sinusoidal profile | | Random field | |
|---|---|---|---|---|
| | $R^2$ | Relative $L_2$ error | $R^2$ | Relative $L_2$ error |
| 0.01 | 0.5987 | 0.6335 | 0.3227 | 0.8230 |
| 0.05 | 0.6446 | 0.5962 | 0.3650 | 0.7969 |

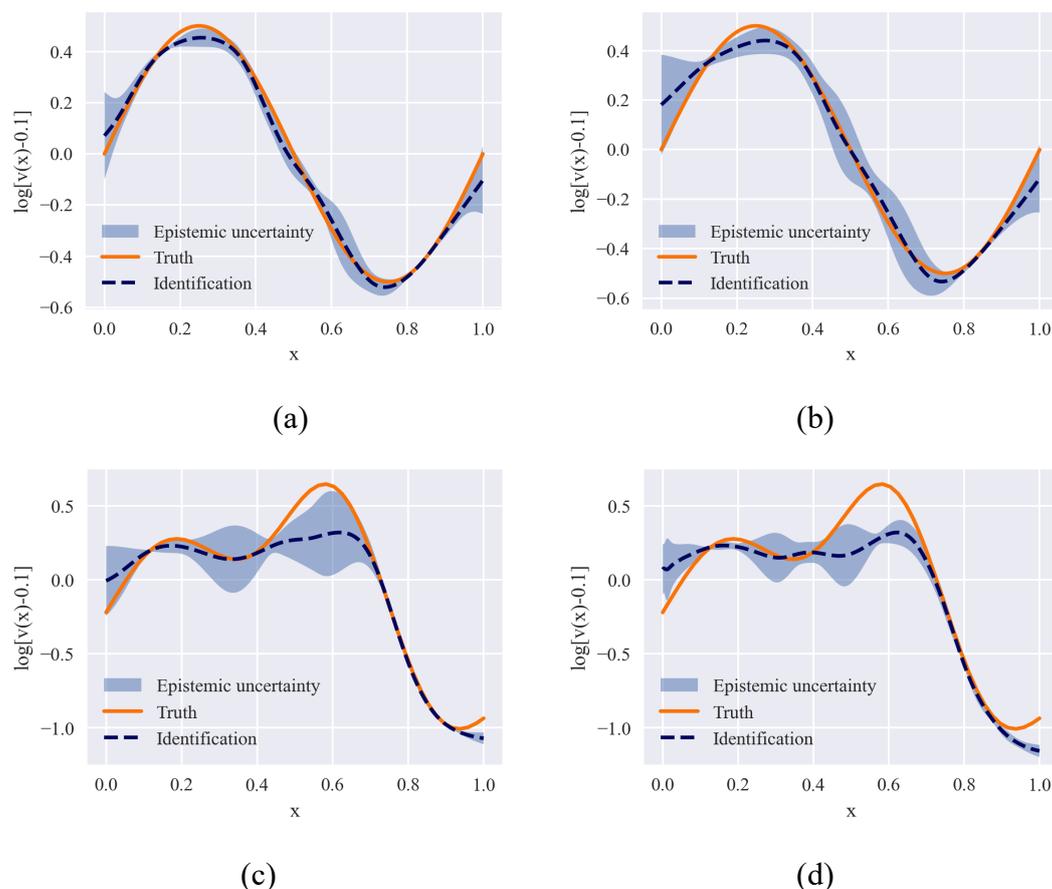

(a)       (b)

(c)       (d)

**Figure 6.** Inverse results of E-PINNs trained with AT. (a) Sinusoidal profile, noise 0.01; (b) Sinusoidal profile, noise 0.05; (c) Random field, noise 0.01; (d) Random field, noise 0.05.

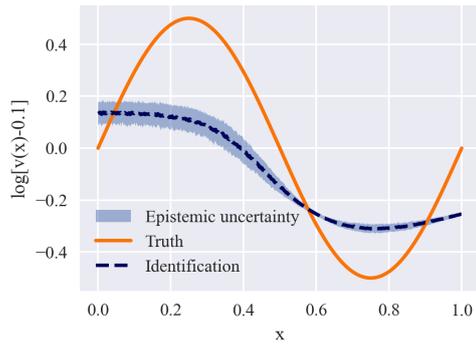
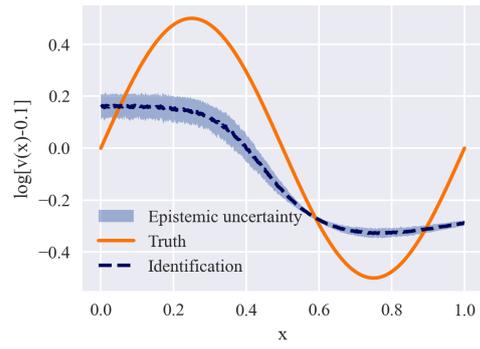

(a)                  (b)

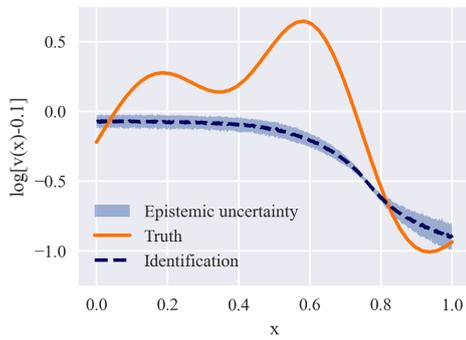
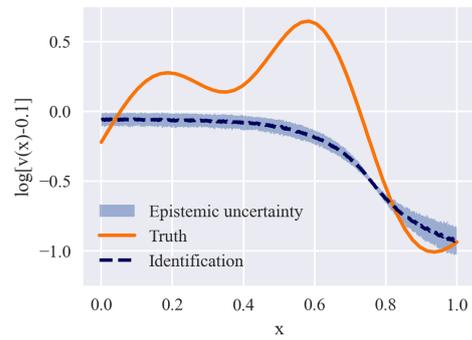

(c)                  (d)

**Figure 7.** Inverse results of MC-dropout. (a) Sinusoidal profile, noise 0.01; (b) Sinusoidal profile, noise 0.05; (c) Random field, noise 0.01; (d) Random field, noise 0.05.

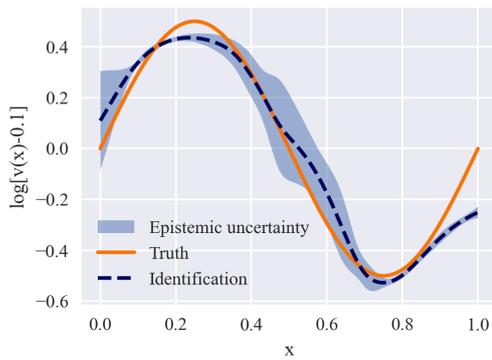
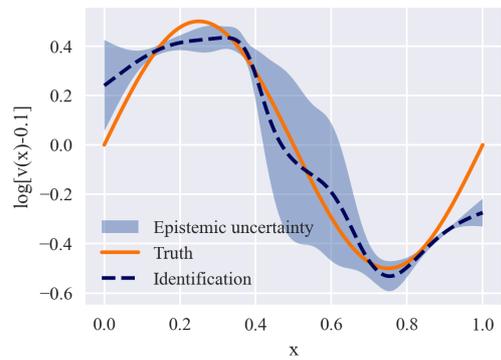

(a)                  (b)

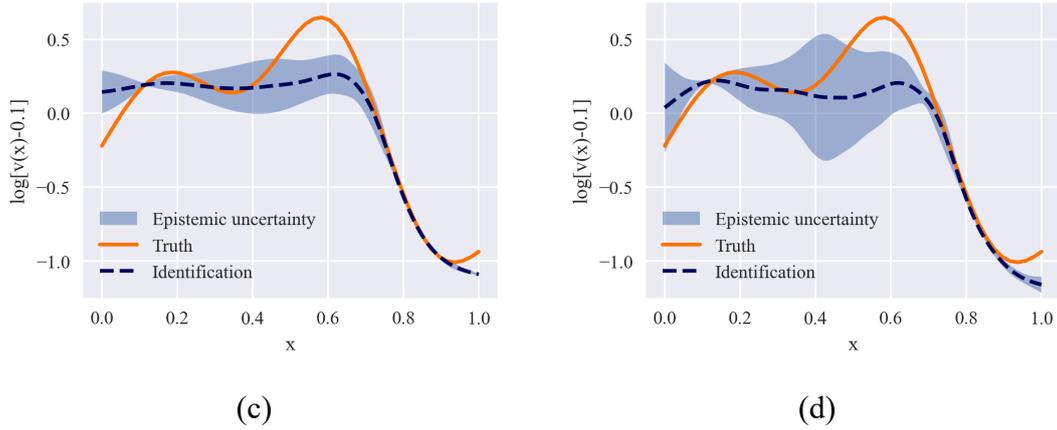

(c)                          (d)

**Figure 8.** Inverse results of E-PINNs without AT. (a) Sinusoidal profile, noise 0.01; (b) Sinusoidal profile, noise 0.05; (c) Random field, noise 0.01; (d) Random field, noise 0.05.

4.2.2 *Adaptive active sampling procedure*

In this section, an adaptive AS strategy is used to demonstrate the use of high-quality uncertainty estimates obtained by E-PINNs and to improve the accuracy of inversion. The cases with the noise of 0.01 are taken as an example to show the effectiveness of this strategy. The implementation of AS is based on the E-PINN trained in the previous section. Since the members of E-PINNs are well-trained, members could reach new convergence faster during the procedure of AS. Experiments have found that 10,000 training steps are enough to ensure convergence in experiments of this case, while other hyperparameters of Adam remain unchanged.

Figure 9 and Figure 10 show the iterative process of gradually increasing the observed sample. Here the scaling factor $\alpha$ is set to 0.01 to obtain sufficiently accurate results. The cases of the sinusoidal profile and the random field stop at steps 8 and 10, respectively. The final inverse result of the sinusoidal profile reaches a very high $R^2$ of 0.9996 and a relative $L_2$ error of 0.0188. The latter reaches an $R^2$ of 0.9995 and a relative $L_2$ error of 0.0218. It is clear that after the first few iterations, the accuracy of inverse results of two cases has reached a fairly high level, the credibility interval becomes very small, and the subsequent iterations have not resulted in much improvement. Considering that the cost of measuring materials may be expensive, it seems reasonable to consider stopping AS in advance if the cost results in only a marginal improvement. For the sinusoidal profile, stopping early seems to be more important since the $R^2$ has reached 0.99 before AS. If $\alpha$ is widened to 0.02, the adaptive AS stops at step 4 for both two cases. While $\alpha$ is set to 0.05, the sampling strategy stops at step 2 for both two cases. According to the results, $\alpha$ is suggested to be set to 0.02 since the sampling cost is greatly reduced, while the R-square of 0.9965 and 0.9976 is still obtained, respectively. Figure 11 and Figure 12 show the evolution trend of the accuracy along with the AS process. From the point of view of $R^2$, it is also appropriate to stop after step 4 since the accuracy of the model has reached a high level, which illustrates the necessity of relaxing the stop criterion. As can be seen from the relative $L_2$ error diagram, although there is a convergence trend, it is difficult to give a reasonable convergence condition in a few iteration steps. Thus, it is reasonable to limit the maximum uncertainty of the results to a certain size as the stopping criterion, which can achieve great accuracy improvement in a few sampling points. In a word, the numerical cases indicate that this strategy could greatly improve the accuracy of inversion in the case of large uncertainty estimates. However, it is necessary to trade-off between accuracy and efficiency.

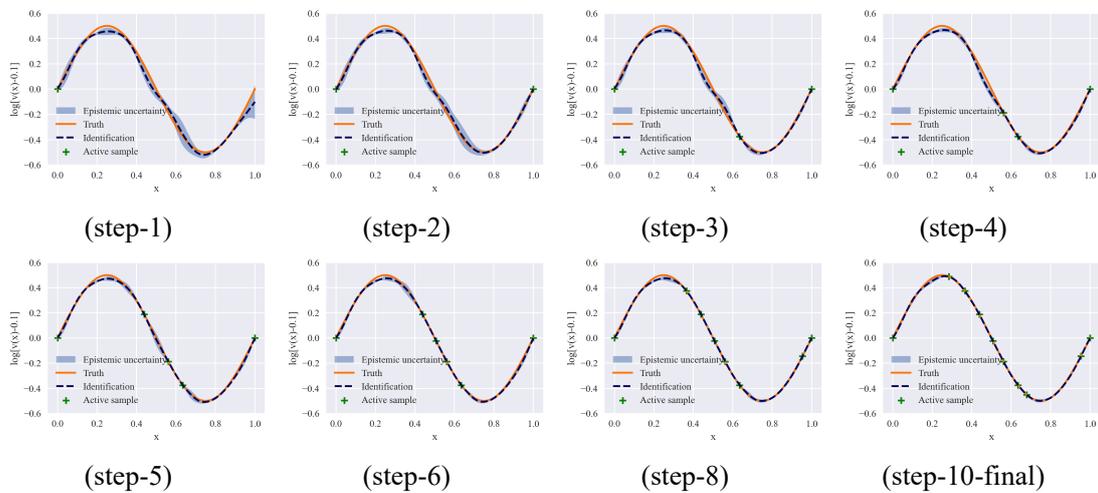

**Figure 9.** The AS for the sinusoidal profile.

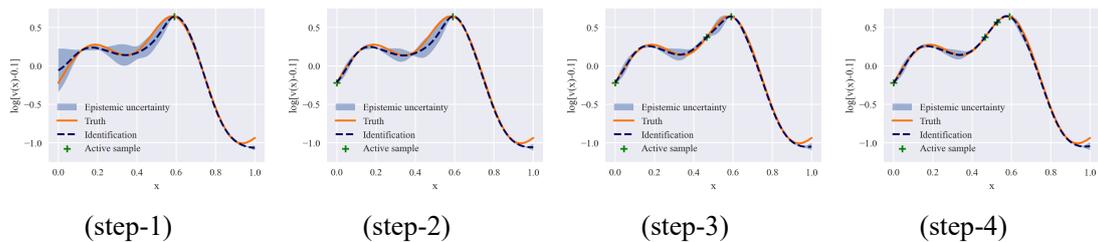

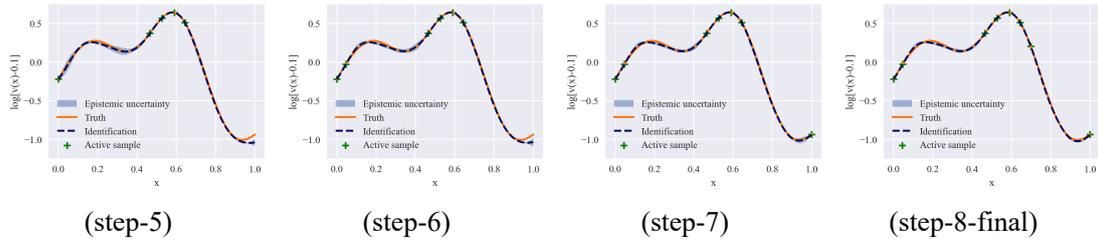

(step-5)      (step-6)      (step-7)      (step-8-final)

**Figure 10.** The AS for the random field.

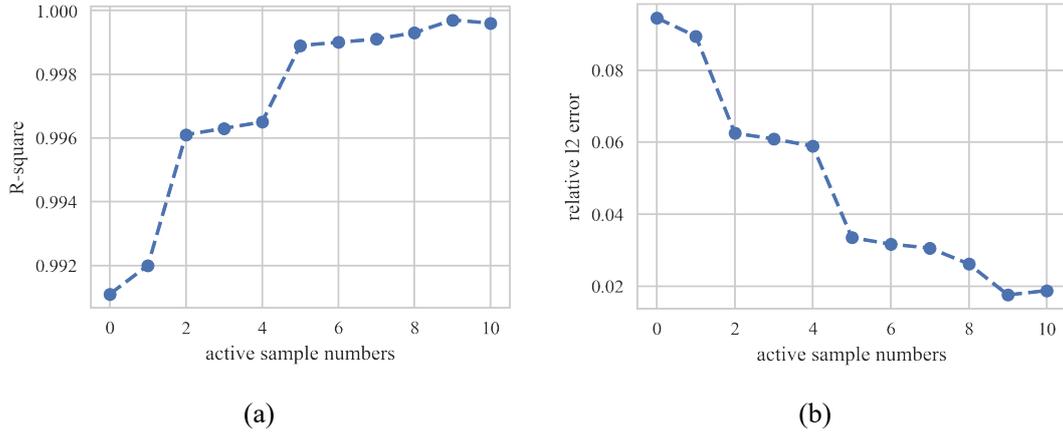

(a)      (b)

**Figure 11.** The AS process of the sinusoidal profile.

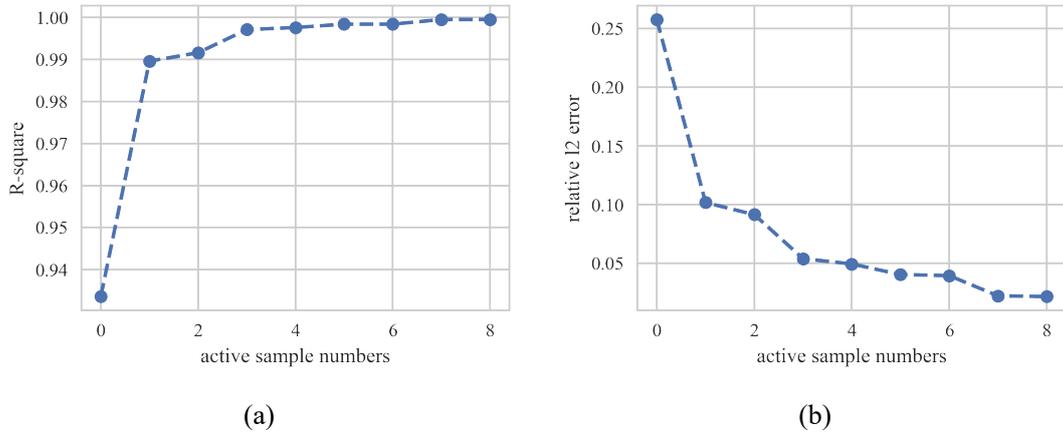

(a)      (b)

**Figure 12.** The AS process of the random field.

## 5. Conclusion

This paper presents a novel method named E-PINN to solve inverse problems whose unknown quantities are fields. The unknown quantity of this kind of inverse problem is propagated by PDEs. The E-PINN employs PINNs as basic models and uses ensemble statistics of several basic models to provide uncertainty quantifications for the solution. This method not only avoids the repeated forward problem evaluations in the classical inverse method (*e.g.,* MCMC) but also robustly identifies the QoI from the noisy data. From a practical viewpoint, the E-PINN offers several advantages. It is straightforward to apply and no sophisticated tools are necessary. There is no need to make too many complex prior assumptions (*e.g.,* likelihood function) about the problem beforehand. Several numerical examples are presented and the results are compared with MC-dropout. The results demonstrate the effectiveness of E-PINNs. Specifically, the main findings of this study can be summarized as follows:

- This study explored the E-PINN for solving inverse problems and found that a small number of model $M = 5$ is enough to obtain high-quality uncertainty estimates. When the number of models is increased to 10, no significant improvement in accuracy and uncertainty estimate is observed in the studied numerical examples. The E-PINN aims at quantifying epistemic uncertainty, and it bypasses the challenge of Bayesian inference and uses ensemble statistics to make uncertainty estimates. No prior assumptions are made for a quantity of interest (QoI) that might limit the expressiveness of the model. Only a few hyperparameters need to be tuned. These features of the E-PINN make it straightforward to be applied. Furthermore, each member of E-PINNs is trained independently, so that the training process of this method can be easily parallelized.

- The FGM is introduced to training the surrogate of the response field (observation data). The results show that AT improves the robustness of the model and improves the accuracy of model inversion to varying degrees. AT is a regularization technique, which might not play an important role in all problems, but it deserves to be employed to solve inverse problems stably. According to the results, AT can also improve the uncertainty estimate quality of E-PINNs. It might be because it avoids overfitting to a certain extent so that better model diversity can be obtained in the ensemble.

- As for the material field inverse problem, pretraining is used to accelerate the convergence of the training process. Based on the uncertainty estimates of E-PINNs, the accuracy of material field inversion is improved by using AS, and an adaptive stopping criterion for AS is proposed. Sampling points are determined according to the maximum uncertainty of the inversion results, and sampling stops when the maximum uncertainty is lower than the predefined threshold. However, it is worth noting that the stop criterion for AS requires a trade-off between accuracy and cost. As for the problem of high sampling and computation cost, the convergence condition should be relaxed to reduce the cost.

## Acknowledgments

This work has been supported by Project of the National Natural Science Foundation of China (11972155), Peacock Program for Overseas High-Level Talents Introduction of Shenzhen City (KQTD20200820113110016) and Project supported by Provincial Natural Science Foundation of Hunan(2020JJ4945).

## Declaration of interests

The authors declared that they have no conflicts of interest with this work. We declare that we do not have any commercial or associative interest that represents a conflict of interest in connection with this work.